\documentclass[11pt]{amsart}
\usepackage{amsbsy}
\usepackage{amsmath}
\usepackage{amsbsy}
\usepackage{graphicx}
\usepackage{enumerate}
\usepackage[cp1250]{inputenc}
\usepackage{indentfirst}
\usepackage{latexsym}
\usepackage{amssymb}
\usepackage{amsfonts}
\usepackage{color}
\usepackage{dsfont}
\usepackage{mathtools}
\usepackage{amsthm}
\usepackage{tikz-cd}

\usepackage{tikz}
\usetikzlibrary{arrows}
\usepackage{xcolor}
\usepackage{xypic}
\usepackage{cite}
\usepackage{float}
\usepackage{filecontents}
\usepackage{graphicx}
\usepackage{epstopdf}
\usepackage{epsfig}

\definecolor{darkgreen}{rgb}{0,.7,.3}

\theoremstyle{definition}
\newtheorem{definition}{Definition}[section]
\newtheorem{remark}[definition]{Remark}

\theoremstyle{plain}
\newtheorem{lemma}[definition]{Lemma}
\newtheorem{theorem}[definition] {Theorem}

\newtheorem{corollary}[definition]{Corollary}

\newtheorem{assumption}[definition]{Assumption}

\title[Finite non-parabolic subgroups]
{Finite non-parabolic subgroups of relatively hyperbolic groups}

\author{Oleg Bogopolski}
\address{Institute of Mathematics, University of Szczecin, ul. Wielkopolska 15, 70-451 Szczecin, Poland}
\email{oleg.bogopolskiy@usz.edu.pl}

\begin{document}


\begin{abstract} 
Let $G$ be a group that is relatively hyperbolic with respect to a collection of subgroups $\{H_{\lambda}\}_{\lambda\in \Lambda}$. Suppose that $G$ is given by a finite relative presentation $\mathcal{P}$ with respect to this collection. We give an upper bound on the orders of finite non-parabolic subgroups of $G$ in terms of some fundamental constants associated with $\mathcal{P}$. This upper bound is computable if $G$ is finitely generated and the word problem in each $H_{\lambda}$, $\lambda\in \Lambda$, is decidable. 
\end{abstract}

\maketitle


\section{Introduction}

Relatively hyperbolic groups were first introduced by Gromov~\cite{Gromov} and then studied by Bowditch~\cite{Bowditch}, Farb~\cite{Farb}, Osin~\cite{Osin_0}, and many other authors from different points of view. Examples of these groups include fundamental groups of finite-volume noncompact Riemannian manifolds of pinched negative curvature, geometrically finite Kleinian groups, word hyperbolic groups, and small cancellation quotients of free products. In this paper, we refer to the Osin manuscript ~\cite{Osin_0}, especially because of the technique developed there.  
 
It is well known that if $G$ is a hyperbolic group and $X$ is a finite generating set of $G$, then the order of any finite subgroup of $G$ is bounded from above by $\bigl(2|X|)^{4\delta+2}$, where $\delta$ is the hyperbolicity constant of the Cayley graph $\Gamma(G,X)$. This easily follows from the result of Bogopolski and Gerasimov~\cite{Bog_Ger}, see Subsection~2.2.
In this paper, we give an upper bound on the orders of finite non-parabolic subgroups of relatively hyperbolic groups. 
\begin{definition} Let $G$ be a group that is relatively hyperbolic with respect to a collection of subgroups
$\{H_{\lambda}\}_{\lambda\in \Lambda}$. A subgroup $F$ of $G$ is called {\it non-parabolic} if there is no $\lambda\in \Lambda$ such that $F$ is conjugate to a subgroup of $H_{\lambda}$. 
\end{definition}

Other notions used in the following theorem are given in Section 2. 

\begin{theorem}\label{finite_subgroups_rel_hyperb_groups}
Let $G$ be a group that is relatively hyperbolic with respect to a collection of subgroups $\{H_{\lambda}\}_{\lambda\in \Lambda}$. Suppose that $G$ is given by a finite relative presentation
$\mathcal{P}=\langle X\sqcup \mathcal{H}\,|\, \mathcal{R}\sqcup \mathcal{Q}\rangle $. Then
the order of any finite non-parabolic subgroup of $G$ is bounded from above by 
$$
\bigl((2|X\cup \Omega|)^{(2M\mathcal{C}+1)(4\delta+2)^2}\bigr) !,
$$
where 
\begin{enumerate}
\item[\tiny{\textcolor{gray}{$\blacktriangleright$}}] 
$\mathcal{C}\geqslant 1$ is the integer relative isoperimetric constant associated with $\mathcal{P}$,

\item[\tiny{\textcolor{gray}{$\blacktriangleright$}}]  $\delta\geqslant 1$ is the integer
hyperbolicity constant of the Cayley graph $\Gamma(G,X\sqcup \mathcal{H})$,

\item[\tiny{\textcolor{gray}{$\blacktriangleright$}}] $\Omega$ is the set of all $\mathcal{H}$-letters in all relators from $\mathcal{R}$.



\item[\tiny{\textcolor{gray}{$\blacktriangleright$}}] $M$ is the maximum of lengths of relators from $\mathcal{R}$. 

\end{enumerate}

This upper bound is computable if Assumption~\ref{assumprion_rel_hyp} is satisfied.

If all $H_{\lambda}$'s are torsion-free, the factorial in the above estimate can be removed. 
\end{theorem}


\begin{remark}
To our best knowledge, no upper bound was known even in the case of finite cyclic non-parabolic subgroups.
However, it was known that the number of conjugacy classes of finite order non-parabolic elements of $G$ is finite (see~\cite[Theorem~4.2]{Osin_0}). 
\end{remark}

Recall that the {\it order problem for elements} in a group $G$ is to calculate the order of a given element of $G$.

Suppose that $G$ is a finitely generated group relatively hyperbolic with respect to a collection of subgroups $H_1,\dots,H_m$, and that the word problem is solvable in each~$H_i$.
Theorem 5.15 in~\cite{Osin_0} says that in this case the order problem for non-parabolic elements in $G$ is solvable.
However, the proof of this theorem, given in~\cite{Osin_0}, 
establishes only the {\it existence} of an algorithm for each {\it concrete} $G$.
In contrast, the algorithm in the proof of the following corollary is a computing algorithm in the class of all such $G$'s, that is, our algorithm is constructive and universal.

\begin{corollary}\label{corollary}
There exists an algorithm solving the order problem for non-parab\-olic elements in the class of all finitely generated relatively hyperbolic groups given as in Assumption \ref{assumprion_rel_hyp}. More precisely: 
The input of this algorithm is

\begin{enumerate}
\item[--] a finite set $X$ (a finite relative generating set of $G$).

\item[--] finite sets $Y_1,\dots ,Y_m$ (finite generating sets of the peripheral subgroups\break $H_1,\dots,H_m$, respectively),

\item[--] algorithms $\mathcal{A}_1,\dots,\mathcal{A}_m$ solving the word problem in $H_1,\dots,H_m$, respectively, 

\item[--] a finite set $\mathcal{R}$ of relative relations of $G$ written in the alphabet $Z\sqcup Z^{-1}$, where $Z=X\sqcup \overset{m}{\underset{i=1}{\sqcup}} Y_i$,

\item[--] a word $W$ in the alphabet $Z\sqcup Z^{-1}$ representing a non-parabolic element\break of $G$.
\end{enumerate}

The output is the order of this element in $G$.  
\end{corollary}

We use this corollary in~\cite{Bogo_Bier} as part of an algorithm for solving exponential equations in relatively hyperbolic groups. Note that it is decidable whether a given word in the alphabet $Z\sqcup Z^{-1}$ represents a parabolic element of $G$ or not (see Theorem~5.6 in~\cite{Osin_0}).

The structure of the paper is as follows.
All basic notations and facts are given in Section 2. In Section 3, we introduce notions that are crucial for
the proof of Theorem~\ref{finite_subgroups_rel_hyperb_groups}. Section~4 contains an outline of the proof and lemmas on decreasing the norms of certain subsets of finite subgroups of $G$, and Section 5 contains the proofs of the theorem and corollary.

\section{Preliminaries}

\subsection{General notations.} 
Given a set $\mathcal{Z}$, we denote by $\mathcal{Z}^{\ast}$ the free monoid generated by $\mathcal{Z}$. 
The length of a word $W\in \mathcal{Z}^{\ast} $ is denoted by $||W||$. 

Let $G$ be a group generated by a subset $Y$. For $g\in G$, let $|g|_Y$ be the length of any shortest word in the alphabet $Y\sqcup Y^{-1}$ representing $g$. 
The right Cayley graph of $G$ with respect to $Y$ is denoted by $\Gamma(G,Y)$. 
By a path $p$ in the Cayley graph we mean a combinatorial path. The initial and terminal vertices of $p$ are denoted by $p_{-}$ and $p_{+}$, respectively. The path inverse to $p$ is denoted by $\overline{p}$.
The length of $p$ is denoted by $\ell(p)$. The {\it label} of $p$  (which is a word in the free monoid  $(Y\sqcup Y^{-1})^{\ast}$)
is denoted by ${\mathbf{ Lab}}(p)$.
The canonical image of ${\mathbf{Lab}}(p)$ in $G$ is denoted by ${\mathbf{Lab}}_G(p)$.

\subsection{Hyperbolic spaces and finite subgroups of hyperbolic groups}
In this paper, we use the definition of a $\delta$-hyperbolic space via the Rips condition. That is, a metric space $S$ is {\it $\delta$-hyperbolic} if it is geodesic and for any geodesic triangle $\Delta$ in $S$, every side of $\Delta$ is contained in the union of the $\delta$-neighborhoods of the other two sides. A finitely generated group $G$ is called $\delta$-{\it hyperbolic} with respect to a finite generating set $Y$ if the Cayley graph $\Gamma(G,Y)$ is $\delta$-hyperbolic. 

There is an equivalent definition of hyperbolicity in terms of $\delta_1$-thin geodesic triangles, see~\cite{BH}. In~\cite{Bog_Ger}, we proved that if $G$ is a finitely generated group that is $\delta_1$-thin with respect to a finite generating set $Y$, then any finite subgroup of $G$ 
can be conjugated into the ball of radius $2\delta_1+1$ and center $1$ in $\Gamma(G,Y)$. Note that the proof does not use the finiteness of $Y$.
Since $\delta_1$ can be taken equal to $2\delta$ (see the proof of Proposition 1.17 in~\cite[Ch. III.H]{BH}), we can reformulate this result as follows:

\begin{theorem} {\rm (see~\cite{Bog_Ger})}\label{Bog_Ger}
If $G$ is a group and $Y$ is a (not necessarily finite) generating set of $G$ such that the Cayley graph $\Gamma(G,Y)$ is $\delta$-hyperbolic, then any finite subgroup of $G$ can be conjugated into the ball of radius $4\delta+1$ and center $1$ in $\Gamma(G,Y)$.
\end{theorem}

An alternative proof of this theorem is given in~\cite{BH}.
In the present paper, we will use this result in the situation where $Y$ is infinite. 

\subsection{Relative presentations of groups, relative areas, and relative Dehn functions}
Let $G$ be a group and $\{H_{\lambda}\}_{\lambda\in \Lambda}$ a collection of subgroups of $G$.
A subset $X$ of $G$ is called a {\it relative generating set of $G$ with respect to}
$\{H_{\lambda}\}_{\lambda\in \Lambda}$ if $G$ is generated by $X$ together with the union of all $H_{\lambda}$'s.
We set 
$$
\mathcal{X}=X\sqcup X^{-1}, \hspace*{4mm} \mathcal{H}=\underset{\lambda\in \Lambda}{\sqcup} (H_{\lambda}\setminus\{1\}),\hspace*{3mm}{\text{\rm and}}\hspace*{2mm} \mathcal{Q}=\underset{\lambda\in \Lambda}{\sqcup} \mathcal{Q}_{\lambda},
$$ where $\mathcal{Q}_{\lambda}$ is the set of all words of length at most 3 in the alphabet $H_{\lambda}\setminus\{1\}$ that represent 1 in the group $H_{\lambda}$ (think about the Cayley multiplication table for $H_{\lambda}$). 

\begin{remark}
Note that in the definition of Osin (see~\cite[Section~2.1]{Osin_0}),
each $\mathcal{Q}_{\lambda}$ consists of all
finite words in the alphabet $H_{\lambda}\setminus\{1\}$ that represent 1 in the group $H_{\lambda}$. However, these two approaches are equivalent.
\end{remark}

A subset $\mathcal{R}\subseteq (\mathcal{X}\sqcup \mathcal{H})^{\ast}$ 
is called a {\it relative set of defining relations} in $G$, if $G$ is canonically isomorphic to the quotient of the free group $F(X\sqcup\mathcal{H})$ by the normal closure of the set $\mathcal{R}\cup \mathcal{Q}$. In this case 
\begin{equation}\label{present}
\mathcal{P}=\langle X\sqcup \mathcal{H}\,|\, \mathcal{R}\sqcup \mathcal{Q}\rangle
\end{equation}
is called a {\it relative presentation} of $G$ with respect to  $\{H_{\lambda}\}_{\lambda\in \Lambda}$. This relative presentation is called {\it finite} if $X$ and $\mathcal{R}$ are finite.
The relative presentation~\eqref{present} is called {\it reduced} if for any $R\in \mathcal{R}$ the cyclic word $R$ does not contain subwords of the form $h_1h_2$, where $h_1,h_2\in H_{\lambda}$ for some $\lambda\in \Lambda$. 

Let $q$ be a cycle in the Cayley graph $\Gamma(G, X\sqcup \mathcal{H})$.
The {\it combinatorial area} (or simply the {\it area}) of $q$ with respect to presentation~\eqref{present}, $Area(q)$, is defined as the minimal number of 2-cells (i.e., of $\mathcal{R}$-cells and $\mathcal{Q}$-cells) in $\Delta$ among all van Kampen diagrams $\Delta$ over presentation~\eqref{present} with the boundary label ${\mathbf{Lab}}(\partial \Delta)\equiv{\mathbf{Lab}}(q)$ (see details in~\cite{Osin_0}; basic notions regarding diagrams can be found in~\cite{Olshanski}). The {\it relative area}, $Area^{rel}(q)$ can be defined similarly if only $\mathcal{R}$-cells are counted. Note that this terminology corresponds to Definition~2.26 from the manuscript of Osin \cite{Osin_0} and differs from the terminology of Dahmani~\cite{Dahmani}, who calls $Area(q)$ the relative area of~$q$.  
These two types of areas are well related: 
$$
Area^{rel}(q)\leqslant Area(q)\leqslant (M+1) Area^{rel}(q)+2\ell(q),
$$
where 
$$
M=\underset{R\in \mathcal{R}}{\max}||R||.
$$
The first inequality is obvious, while the second is a kind of folklore (see 
Lemma~2.52 from~\cite{Osin_0}, which is not the same, but similar). 

\medskip

The {\it Dehn function} for the relative presentation~\eqref{present} is a partial function $D_{\mathcal{P}}$ from $\mathbb{N}$ to $\mathbb{N}$ such that
$$
D_{\mathcal{P}}(n)=\max\{ Area(q)\,|\, q\hspace*{2mm} {\text{\rm is a cycle in}}\hspace*{2mm} \Gamma(G, X\sqcup \mathcal{H})\hspace*{2mm} {\text{\rm of length}}\hspace*{2mm} \ell(q)\leqslant n\}
$$
if this maximum exists; otherwise $D_{\mathcal{P}}(n)$ is not defined.

The {\it relative Dehn function} $D^{rel}_{\mathcal{P}}(n)$ for the relative presentation~\eqref{present} can be defined analogously if we replace $Area(q)$ by $Area^{rel}(q)$.

\subsection{Isolated components}

The following definitions are borrowed from~\cite{Osin_0}. 
Let $q$ be a path (cycic or not) in the Cayley graph $\Gamma(G,X\sqcup \mathcal{H})$. A nontrivial subpath $p$ of $q$
is called an {\it $H_{\lambda}$-subpath} if the label of $p$ is a word in the alphabet $H_{\lambda}$.
An $H_{\lambda}$-subpath $p$ of $q$ is called an {\it $H_{\lambda}$-component} of $q$ if $p$ is not contained in a longer subpath of $q$ with this property. By an {\it $\mathcal{H}$-component} of $q$ we mean an $H_{\lambda}$-component of $q$ for some $\lambda\in \Lambda$. Two $H_{\lambda}$-components $p_1,p_2$ of $q$ are called {\it connected} if there exists a path $\gamma$ in $\Gamma(G,X\sqcup \mathcal{H})$ that connects some vertex of $p_1$ to some vertex of $p_2$ and $\mathbf{Lab}(\gamma)$ is a word from $(H_{\lambda}\setminus \{1\})^{\ast}$. 
Every such $\gamma$
is called a {\it connector} for~$q$.
Note that we can always assume that $\gamma$ has length at most~1, as every element of $H_{\lambda}\setminus \{1\}$
is included in the set of generators of $G$ and the empty word represents 1 by agreement. An $H_{\lambda}$-component $p$ of $q$ is called {\it isolated} if it is not connected to any other $H_{\lambda}$-component of~$q$.

The following definition is slightly different from Definition~2.25 in~\cite{Osin_0}. 

\begin{definition}\label{omega} 
For every $\lambda\in \Lambda$, we denote by $\Omega_{\lambda}$ the set of all $H_{\lambda}$-letters in all words from $\mathcal{R}$. We also put
$$\Omega=\underset{\lambda\in \Lambda}{\cup} \Omega_{\lambda}.
$$
\end{definition}

\noindent
Note that $\Omega$ is finite if $\mathcal{R}$ is finite. Different letters from $\Omega$ can represent the same element of $G$. For simplicity, we will denote this element by any of these letters.

The following lemma is a variation of Lemma~2.27 from~\cite{Osin_0}.

\begin{lemma}\label{omega_length}
Suppose that a group $G$ is given by the finite relative presentation~\eqref{present} 
with respect to a collection of subgroups $\{H_{\lambda}\}_{\lambda\in \Lambda}$. Let $q$ be a cycle in $\Gamma(G,X\sqcup \mathcal{H})$ and let $p_1,\dots,p_k$ be the set of all isolated $\mathcal{H}$-components of $q$. Suppose that $p_i$ is an $H_{\lambda(i)}$-component of $q$, where $i=1,\dots,k$. Then
$\mathbf{Lab}_G(p_i)\in \langle\Omega_{\lambda(i)}\rangle$. Moreover, the lengths of the elements $\mathbf{Lab}_G(p_i)$ with respect to the generating sets $\Omega_{\lambda(i)}$ of the subgroups $\langle \Omega_{\lambda(i)}\rangle$ satisfy the inequality
$$
\overset{k}{\underset{i=1}{\sum}}|\mathbf{Lab}_G(p_i)|_{\Omega_{\lambda(i)}}\leqslant M\cdot Area^{rel}(q),
$$
where $M=\underset{R\in \mathcal{R}}{\max}||R||$ and $Area^{rel}(q)$ is the relative area of $q$.
\end{lemma}

Comparing our Lemma~\ref{omega_length} 
and Lemma 2.27 from~\cite{Osin_0},
we note three differences.
We use another definition of $\Omega_{\lambda}$ and (therefore) we do not assume that the relative presentation~\eqref{present} is reduced.
Moreover, we do not assume that the components $p_1,\dots,p_k$ are $H_{\lambda}$-components for the same $\lambda$. However, the proof of this variation is almost the same.  

\subsection{Relatively hyperbolic groups}
Let $G$ be a group and $\{H_{\lambda}\}_{\lambda\in \Lambda}$ be a system of its subgroups. Recall that $G$ is called a {\it relatively hyperbolic group} with respect to $\{H_{\lambda}\}_{\lambda\in \Lambda}$ if it has a finite relative presentation $\mathcal{P}=\langle X\sqcup \mathcal{H}\,|\, \mathcal{R}\sqcup \mathcal{Q}\rangle$
such that its relative Dehn function $D^{rel}_{\mathcal{P}}$  is defined on the whole $\mathbb{N}$ and bounded from above by a linear function (see~\cite[Definition 2.35]{Osin_0}), i.e.,
there exists a constant $\mathcal{C}>0$ such that
$
D^{rel}_{\mathcal{P}}(n)\leqslant \mathcal{C}n
$
for all natural $n$. Any such constant is called a {\it relative isoperimetric constant} for~$\mathcal{P}$.

It is known (see~\cite[Theorem 2.53 and Corollary 2.54]{Osin_0}) that if $G$ is a group with a finite relative presentation $\mathcal{P}=\langle X\sqcup \mathcal{H}\,|\, \mathcal{R}\sqcup \mathcal{Q}\rangle$ and the relative Dehn function $D^{rel}_{\mathcal{P}}$ is defined on the whole $\mathbb{N}$, then the following conditions are equivalent:

\begin{enumerate}
\item[(1)] The function $D^{rel}_{\mathcal{P}}$ is bounded from above by a linear function.

\item[(2)] The function $D_{\mathcal{P}}$ is bounded from above by a linear function.

\item[(3)] The Cayley graph $\Gamma(G,X\sqcup \mathcal{H})$
is $\delta$-hyperbolic for some $\delta\geqslant 0$.
\end{enumerate}

Thus, if
$G$ is hyperbolic relative to $\{H_{\lambda}\}_{\lambda\in \Lambda}$, then
$
D_{\mathcal{P}}(n)\leqslant \mathcal{C}'n
$
for some constant $\mathcal{C}'>0$. Any such $\mathcal{C}'$ is called an {\it isoperimetric constant} for $\mathcal{P}$.
Since
$
D^{rel}_{\mathcal{P}}(n)\leqslant D_{\mathcal{P}}(n),
$
we may always assume that $\mathcal{C}\leqslant \mathcal{C}'$. 


\subsection{Computability of the constants $\mathcal{C}$ and $\delta$}

Many algorithms require knowledge of the constants $\mathcal{C}$ and $\delta$. Fortunately, they can be algorithmically computed under a certain natural assumption about $G$.

\begin{assumption}\label {assumprion_rel_hyp} Let $G$ be a relatively hyperbolic group with respect to a collection of subgroups $\{H_{\lambda}\}_{\lambda\in \Lambda}$.
Suppose that $\Lambda$ is finite and that for any $\lambda\in \Lambda$
the group $H_\lambda$ is generated by a finite set~$Y_{\lambda}$, and the word problem in $H_{\lambda}$ is decidable.
Assume that the following data are given:

\begin{enumerate}
\item[{\rm (1)}] The sets 
$X$ and $\mathcal{R}$ from a finite
relative presentation $\mathcal{P}=\langle X\sqcup\mathcal{H}\,|\, \mathcal{R}\sqcup \mathcal{Q}\rangle$ of $G$ (see Subsection 2.3), where each relator $R\in \mathcal{R}$
is written as a word in the alphabet $Z\sqcup Z^{-1}$, where $Z=X\sqcup \underset{\lambda\in \Lambda}{\sqcup} Y_{\lambda}$.

\item[{\rm (2)}] Algorithms $\mathcal{A}_{\lambda}$ deciding the word problem in $H_{\lambda}$ for words in the alphabet $Y_{\lambda}\sqcup Y_{\lambda}^{-1}$, $\lambda\in \Lambda$.
\end{enumerate}
\end{assumption}

\begin{theorem}\label{theorem_0f_Dahmani} {\rm (see~\cite[Proposition 7.3]{Dahmani})}
Let $G$ be a relatively hyperbolic group given as in 
Assumption~\ref{assumprion_rel_hyp}.
Then there exists an algorithm computing an isoperimetric constant $\mathcal{C}'>0$ 
for the relative presentation $\mathcal{P}$.
\end{theorem}

Knowing the isoperimetric constant $\mathcal{C}'$, one can compute a constant $\delta\geqslant 0$ such that the Cayley graph $\Gamma(G,X\sqcup \mathcal{H})$ is $\delta$-hyperbolic,
see~Theorem 2.9 in Ch.~III.H of~\cite{BH}. 

\begin{corollary}
Let $G$ be a relatively hyperbolic group given as in
Assumption~\ref{assumprion_rel_hyp}.
Then there exists an algorithm computing the hyperbolicity constant $\delta\geqslant 0$ 
of the Cayley graph $\Gamma(G,X\sqcup \mathcal{H})$.
\end{corollary}

The relative isoperimetric constant $\mathcal{C}$ for $\mathcal{P}$ can be taken equal to $\mathcal{C}'$. 
For convenience, we increase the constants $\mathcal{C}$, $\mathcal{C}'$ and $\delta$ so that they become integers. 

\section{Doubly $\Lambda$-reduced words}

Let $G$ be a relatively hyperbolic group with respect to a collection of subgroups $\{H_{\lambda}\}_{\lambda\in \Lambda}$ and let $\mathcal{P}=\langle X\sqcup\mathcal{H}\,|\, \mathcal{R}\sqcup \mathcal{Q}\rangle$
be a finite relative presentation of $G$ with respect to this collection.
Recall that
$$
\mathcal{X}=X\sqcup X^{-1}\hspace*{2mm}{\text{\rm and}}\hspace*{2mm}\mathcal{H}=\underset{\lambda\in \Lambda}{\sqcup} (H_{\lambda}\setminus\{1\}).
$$
In what follows, we work mainly with the free monoid
$(\mathcal{X}\sqcup \mathcal{H})^{\ast}$. 


Below we define $\Lambda$-reduced and doubly $\Lambda$-reduced words in this monoid

\begin{definition}
Let $W$ be a word from the monoid $(\mathcal{X}\sqcup \mathcal{H})^{\ast}$. 
The image of $W$ under the natural map $\varphi:(\mathcal{X}\sqcup \mathcal{H})^{\ast}\rightarrow G$ is denoted by $\overline{W}$. 
For any subset ${\mathcal{S}}$ of $(\mathcal{X}\sqcup \mathcal{H})^{\ast}$,
the image of $\mathcal{S}$ under $\varphi$ is also denoted by $\overline{\mathcal{S}}$.
The empty word is denoted by $\varepsilon$. We have $\overline{\varepsilon}=1$. 
Let $W=w_1\dots w_n$, where $w_1,\dots ,w_n\in \mathcal{X}\sqcup \mathcal{H}$. The {\it formal inverse} to $W$ is the word
$W^{-1}:=w_n^{-1}\dots w_1^{-1}$. 
Recall that the length of $W$ is denoted by $||W||$. 

If $W$ is nonempty, then the first and the last letters of $W$ are denoted by $W_{-}$ and $W_{+}$, respectively.
A nonempty subword $U$ of $W$
is called a {\it syllable} of $W$ if either $U\in \mathcal{X}^{\ast}$,
or $U\in (H_{\lambda}\setminus \{1\})^{\ast}$ for some $\lambda$ in $\Lambda$, and $U$ is not contained in a longer subword of $W$ with this property. In the second case, $U$ is called an {\it $H_{\lambda}$-syllable} (or simply an {\it $\mathcal{H}$-syllable}) of $W$. Now we introduce the main notions. 

\begin{enumerate}
\item[1)] $W$ is called {\it reduced} if it does not contain subwords 
of the form $uv$, where $u,v\in H_{\lambda}\setminus \{1\}$ for some $\lambda\in \Lambda$, i.e., if all $H_{\lambda}$-syllables of $W$ have length 1.
Given a word $U\in (\mathcal{X}\sqcup \mathcal{H})^{\ast}$, there exists a unique reduced word in $(\mathcal{X}\sqcup \mathcal{H})^{\ast}$ that is equal to $U$ in the free product
$F(\mathcal{X})\ast (\underset{\lambda\in \Lambda}{\ast} H_{\lambda})$. We denote it by $U_{\rm red}$. Note that each word in $\mathcal{X}^{\ast}$ is reduced. 

\item[2)] $W$ is called {\it geodesic} if
$W$ is a shortest word among all words representing the element $\overline{W}$, i.e., if $||W||=|\overline{W}|_{X\cup \mathcal{H}}$. In particular, geodesic words are reduced.
 
\item[3)] 
$W$ is called {\it cyclically reduced} if it is reduced and the first and last letters of $W$ do not lie in the same $H_{\lambda}$.
 
Equivalently, $W$ is cyclically reduced if $W^k$ is reduced for any $k\geqslant 1$.

We stress that if $W$ is cyclically reduced, then $W$ is not a letter from $\mathcal{H}$. 

\item[4)] $W$ is called {\it $\Lambda$-reduced} if it is reduced and for any (equivalently, for some) path $p$ in $\Gamma(G, X\sqcup \mathcal{H})$ with the label $W$ all $H_{\lambda}$-components of $p$ are isolated, $\lambda\in \Lambda$. 

Equivalently, $W$ is $\Lambda$-reduced if it is reduced and the following is satisfied. Write $W$ as the product of syllables, $W=u_1u_2\dots u_k$. Then for any two syllables $u_i$ and $u_j$, where $i<j$ and $u_i,u_j\in H_{\lambda}$ for some $\lambda\in\Lambda$, we have $\overline {u_{i+1}\dots u_{j-1}}\notin H_{\lambda}$.

\item[5)] $W$ is called {\it doubly $\Lambda$-reduced} if it is cyclically reduced and $WW$ is $\Lambda$-reduced. 
\end{enumerate}
\end{definition} 

Note that every geodesic word $W$ is $\Lambda$-reduced and that the empty word $\varepsilon$ is doubly $\Lambda$-reduced.  
The following lemma shows that if $f\in G$ is a non-parabolic element of minimal $(X\cup \mathcal{H})$-length in the conjugacy class $f^G$, then every geodesic word $W$ representing $f$ is doubly $\Lambda$-reduced.

\begin{lemma}\label{conj_geodes_doubly_Lambda-reduced}
Let $G$ and $\mathcal{P}$ be as above. Suppose that $W$ is a word from $(\mathcal{X}\sqcup \mathcal{H})^{\ast}\setminus \mathcal{H}$ that is geodesic and not doubly $\Lambda$-reduced. 
Then there exist a word $W_1\in (\mathcal{X}\sqcup \mathcal{H})^{\ast}$ and a nonempty terminal subword $U$ of $W$ such that
$$
\overline{W_1}=\overline{UWU^{-1}}\hspace*{3mm}{\text{\rm and}}\hspace*{3mm}
||W_1||<||W||\hspace*{3mm} {\text{\rm and}}\hspace*{3mm} ||U||<||W||.
$$
If, additionally, $W$ is cyclically reduced, then we can choose the above $U$ so that
$$\overline{U}\in \langle X\cup \Omega\rangle\hspace*{3mm}{\text{\rm and}}\hspace*{3mm} 
|\overline{U}|_{X\cup \Omega}\leqslant (2M\mathcal{C}+1)||W||.
$$
\end{lemma}

{\it Proof.} We write $W=w_1w_2\dots w_n$, where $w_i\in \mathcal{X}\sqcup \mathcal{H}$, $i=1,\dots ,n$.  
We may assume that $W$ is cyclically reduced; otherwise, the lemma is evident. 
Then, by assumption, $W^2$ is not $\Lambda$-reduced.

 Let
$P=e_1e_2\dots e_ne_1'e_2'\dots e_n'$ be any path in $\Gamma(G,X\sqcup \mathcal{H})$ with the label $W^2$, i.e., the labels of the edges $e_i$ and $e_i'$ are $w_i$, $i=1,\dots,n$. 

Since $W$ is cyclically reduced,
the edges $e_n$ and $e_1'$ do not belong to the same component of $P$. Therefore (and since $W$ is geodesic), every component of $P$ is an edge. 
Since $W$ is geodesic and $W^2$ is not $\Lambda$-reduced, $P$ contains two connected $H_{\lambda}$-components $e_i$ and $e_j'$ with $i<n$ and $j>1$. 

Let $e$ be a path in $\Gamma(G,X\sqcup \mathcal{H})$ with label in $(H_{\lambda}\setminus \{1\})^{\ast}$ such that $(e_i)_{+}=e_{-}$ and $(e_j)_{-}'=e_{+}$. We may assume that $e$ is an edge or a trivial path. We show that the statement is valid for $U=w_{i+1}\dots w_n$ and certain $W_1$. 
Consider two cases.

\medskip

{\it Case 1.} Suppose that $i\geqslant j$ (see Fig. 1). 

Consider the paths
$p=e_{i+1}\dots e_ne_1'\dots e'_n\overline{e'_n}\dots \overline{e'_{i+1}}$ and
$q=ee_j'\dots e_i'$. Since $p_{-}=q_{-}$ and $p_{+}=q_{+}$, we have $\mathbf{Lab}_G(q)=\mathbf{Lab}_G(p)$. Thus,
$\overline{W_1}=\overline{UWU^{-1}}$ for $W_1=\mathbf{Lab}(q)$.
Moreover,
$$
||W_1||=\ell (q)\leqslant
i-j+2
\leqslant (n-1)-2+2
<n=||W||.
$$

{\it Case 2.} Suppose that $j>i$ (see Fig. 2).

This case can be considered similarly using the path $q=e\overline{e_{j-1}'}\dots \overline{e_{i+1}'}$ and the word $W_1={\bf Lab}(q)$. 
As above, we have $\overline{W_1}=\overline{UWU^{-1}}$.
Note that $j\neq i+1$ since $W$ is reduced. Then
$$
||W_1||=\ell (q)\leqslant (j-1)-(i+1)+2=j-i\leqslant n-1<n=||W||.
$$

\medskip

\vspace*{2mm}
\hspace*{2mm}
\includegraphics[scale=0.6]{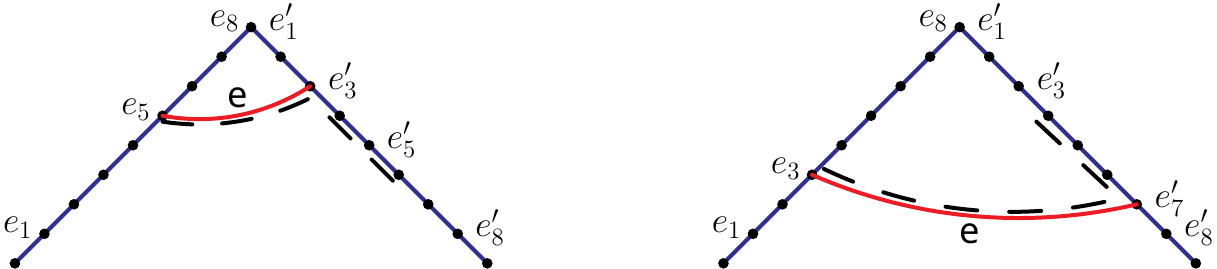}

\vspace*{5mm}

\begin{center}
Fig. 1. Case 1 for $i=5$, $j=3$.\hspace*{18mm}  Fig. 2. Case 2 for $i=3$, $j=7$.
\end{center}

\medskip

Now we prove the second statement of the lemma. 
Passing to another pair of connected components $e_i, e_j'$ of $P$ if necessary, we may assume that all components of the cycle $Q=e_{i+1}\dots e_ne_1'\dots e_{j-1}'\overline{e}$ are isolated edges. 
By Lemma~\ref{omega_length}, the labels in $G$ of all components of $Q$ lie in $\langle \Omega \rangle$ and the sum of $\Omega$-lengths of these labels is at most $M\mathcal{C}\ell(Q)\leqslant M\mathcal{C}(2n-1)$. Hence, 
$$
|\overline{U}|_{X\cup \Omega}=|\overline{w_{i+1}\dots w_n}|_ {X\cup \Omega}\leqslant (n-i)+M\mathcal{C}(2n-1)\leqslant (2M\mathcal{C}+1)n= (2M\mathcal{C}+1)||W||.
$$

\hfill $\Box$

\medskip

The following lemma is similar to Lemma 4.5 in the Osin monograph~\cite{Osin_0}. However,
the assumptions and conclusions in these lemmas are different.
Therefore, and for completeness, we give a proof.

\begin{lemma}\label{estimate_Osin} Let $G$ and $\mathcal{P}$ be as above. Let $f\in G$ be an element of finite order and $W\in (\mathcal{X}\sqcup \mathcal{H})^{\ast}$ be a word representing $f$. Suppose that
$W$ is doubly $\Lambda$-reduced. 
Then $f\in \langle X\cup \Omega\rangle$ and
\begin{equation}\label{First_estimate_Osin}
 |f|_{X\cup \Omega}\leqslant (2M\mathcal{C}+1) ||W||.
\end{equation}
Furthermore, if $U_1,\dots , U_s$ are all $\mathcal{H}$-syllables of $W$, then $\overline{U_i}\in \langle \Omega\rangle$ for $i=1,\dots, s$ and
\begin{equation}\label{Second_estimate_Osin}
\overset{s}{\underset{i=1}{\sum}} |\overline{U_i}|_{\Omega}\leqslant 2M\mathcal{C} ||W||.
\end{equation}
\end{lemma}

{\it Proof.} Inequality~\eqref{First_estimate_Osin}
follows straightforwardly from ~\eqref{Second_estimate_Osin}. Therefore, it suffices to prove~\eqref{Second_estimate_Osin}. Let $n$ be the order of $f$.  
We consider an arbitrary cycle $p$ in $\Gamma(G,X\sqcup \mathcal{H})$
with the label $W^n$. If all components of $p$ are isolated, we put $b=p$. If there exist two connected $H_{\lambda}$-components $s_1$ and $s_2$ of $p$, then we proceed as follows. Let $p=as_1bs_2c$. Passing to another pair of connected components if necessary, we may assume that no component of $b$ is connected to the components $s_1$, $s_2$, and all components of $b$ are isolated. 
Since $\mathbf{Lab}(b)$ is a (cyclic) subword of $W^n$, we have $$\mathbf{Lab}(b)\equiv W_0^kV,$$ where $W_0$ is a cyclic shift of $W$ and $V$ is a cyclic subword of $W$ of length less than $||W||$. We have $k>0$, since
$W$ is doubly $\Lambda$-reduced.

Let $e$ be a path in $\Gamma(G,X\sqcup \mathcal{H})$ of length at most 1 such that $e_{-}=b_{-}$, $e_{+}=b_{+}$ and $\mathbf{Lab}_G(e)\in H_{\lambda}$. We consider the cycle $q=be^{-1}$. Note that all components of $q$ are isolated edges. 
Let $q_1,\dots ,q_s$ be all components of the initial subpath of $b$ with the label $W_0$. By Lemma~\ref{omega_length}, applied to the cycle $q$, we have 
$$
\mathbf{Lab}_G(q_i)\in \langle \Omega\rangle$$ for all $i=1,\dots,s$ and
$$
k\overset{s}{\underset{i=1}{\sum}}|\mathbf{Lab}_G(q_i)|_{\Omega}\leqslant M\cdot Area^{rel}(q)\leqslant M\mathcal{C}\ell(q).
$$
The left side of this formula contains the multiple $k$, since $b$ contains at least $k$ ``copies'' of each $q_i$.
We also have
$$
\ell(q) \leqslant \ell(b)+1= k||W_0||
+||V||+1\leqslant (k+1)||W||.
$$
Therefore
\begin{equation}\label{Third_estimate_Osin}
\overset{s}{\underset{i=1}{\sum}}|\mathbf{Lab}_G(q_i)|_{\Omega}\leqslant M\mathcal{C}\frac{k+1}{k} ||W||\leqslant 2M\mathcal{C}||W||.
\end{equation}
Since $W_0$ is a cyclic shift of $W$, the sequence of elements $
\mathbf{Lab}_G(q_1),\dots , \mathbf{Lab}_G(q_s)$, up to a permutation, coincides with $\overline{U_1},\dots,\overline{U_s}$. Thus, \eqref{Third_estimate_Osin} implies~\eqref{Second_estimate_Osin}. 
\hfill $\Box$

\begin{remark} Suppose that 
$f\in G$ is a non-parabolic element of finite order and $f$ can be represented by a word $W\in (\mathcal{X}\sqcup \mathcal{H})^{\ast}$ that is geodesic and cyclically reduced.
In this case
the assumptions of Lemmas~\ref{conj_geodes_doubly_Lambda-reduced} and~\ref{estimate_Osin} 
are mutually complementary. 
\end{remark}

\section{Decreasing norms with conjugations}

In the following proof we use the
constant
$$
K=(2|X\cup\Omega|)^{2M\mathcal{C}+1}.
$$ 
We may assume that $K\geqslant 1$. Indeed, otherwise $X$ and $\Omega$ are empty, hence $\mathcal{R}$ is empty and
$G$ is the free product of the groups $H_{\lambda}$, $\lambda\in \Lambda$. Then $G$ does not have finite non-parabolic subgroups and we are done. 
In this section we also assume that $\Lambda=\{1,\dots,m\}$ for some natural $m\geqslant 1$. The case of arbitrary $\Lambda$ will be considered in Section 5. 

\begin{definition}
For any nonempty subset  $\mathcal{S}\subseteq (\mathcal{X}\sqcup \mathcal{H})^{\ast}$, the cardinal 
$$
\mathcal{N}(\mathcal{S})=\sup \{||W||: W\in \mathcal{S}\}
$$
is called the {\it norm} of $\mathcal{S}$. Analogously, for any nonempty subset $A\subseteq G$, the cardinal
$$
\mathcal{N}_{X\cup \mathcal{H}}(A)=
\sup \{|a|_{X\cup \mathcal{H}}: a\in A\}
$$
is called the {\it norm} of $A$ with respect to $X\cup \mathcal{H}$.
To simplify the notation, we will write $\mathcal{N}(A)$ if the set $X\cup \mathcal{H}$ is fixed. 
We also set $\mathcal{N}(\emptyset)=-\infty$.
\end{definition}

Note that for any subset $\mathcal{S}\subseteq (\mathcal{X}\sqcup \mathcal{H})^{\ast}$ we have $\mathcal{N}(\mathcal{S})
\geqslant \mathcal{N}(\overline{\mathcal{S}})$ with equality if all words in $\mathcal{S}$ are geodesics.

\medskip

Now we can give an outline of the proof of the main theorem. The following two lemmas are devoted to decreasing the norms of certain subsets of $(\mathcal{X}\sqcup \mathcal{H})^{\ast}$. In Section~5 we deduce from these lemmas that if
$H$ is a finite subgroup of the relatively hyperbolic group $G$ and $A$ is a subset of $H$ with norm greater than 1, then a ``large'' part of some conjugate of $A$ has norm smaller than the norm of~$A$. Iterating, we conclude that a ``large'' part of some conjugate of $H$ lies in some peripheral subgroup $H_i$. After that we estimate the cardinal of this part.  

\begin{definition}
For any subsets $\mathcal{S}\subseteq (\mathcal{X}\sqcup \mathcal{H})^{\ast}$ and $A\subseteq G$, we say that $\mathcal{S}$ {\it bijectively represents} $A$ if the map $\mathcal{S}\rightarrow \overline{\mathcal{S}}$, $W\mapsto \overline{W}$, is a bijection and
$\overline{\mathcal{S}}=A$. 
\end{definition}

\begin{lemma}\label{strongly Lambda-reduced}
Let $G$ and $\mathcal{P}$ be as above, and let $H$ be a finite subgroup of $G$.
Suppose that $\mathcal{S}$ is a subset of $(\mathcal{X}\sqcup \mathcal{H})^{\ast}$ that satisfies the following conditions.

\begin{enumerate}
\item[\rm{(a)}] $\mathcal{S}$ is nonempty and different from $\{\varepsilon\}$.

\item[\rm{(b)}] 
$\mathcal{S}$ bijectively represents a subset of~$H$.  

\item[\rm{(c)}] Every word $W\in \mathcal{S}$ is geodesic and cyclically reduced (in particular, $W\notin \mathcal{H}$). 
\end{enumerate}

\medskip

\noindent
Then there exist a word $U\in (\mathcal{X}\sqcup \mathcal{H})^{\ast}$ and a subset $\mathcal{S}_1\subseteq (\mathcal{X}\sqcup \mathcal{H})^{\ast}$ such that

$$
\begin{array}{ll}
||U||< \mathcal{N}(\mathcal{S}), &
\hspace*{5mm}
\mathcal{N}(\mathcal{S}_1)<\mathcal{N}(\mathcal{S}),\vspace*{2mm}\\
\overline{\mathcal{S}_1}\subseteq \overline{U\mathcal{S}U^{-1}},& \hspace*{5mm}
|\mathcal{S}_1|\geqslant \displaystyle{\frac{1}{K^{\mathcal{N}(\mathcal{S})}}}\,\,|\mathcal{S}|-1.
\end{array}
$$
\end{lemma}

\medskip

{\it Proof.}  
We may assume that $|\mathcal{S}|> K^{\mathcal{N}(\mathcal{S})}$, otherwise the conclusion of the lemma is satisfied for $U=\varepsilon$ and $\mathcal{S}_1=\emptyset$.

If some $W\in \mathcal{S}$ is doubly $\Lambda$-reduced,
then, by Lemma~\ref{estimate_Osin}, 
$$\overline{W}\in \langle X\cup \Omega\rangle\hspace*{3mm}{\text{\rm and}}\hspace*{3mm} |\overline{W}|_{X\cup\Omega}\leqslant (2M\mathcal{C}+1)||W||\leqslant (2M\mathcal{C}+1)\mathcal{N}(\mathcal{S}).
$$
Then the number of such $ \overline{W}$ is at most $K^{\mathcal{N}(\mathcal{S})}$. 
Since the representation of $\mathcal{S}$ in $H$ is injective,
the number of such $W$ is also
at most $K^{\mathcal{N}(\mathcal{S})}$.
Therefore, the set 
$$
\mathcal{A}=
\{W\in \mathcal{S}\,|\, W\hspace*{2mm}{\text{\rm is not doubly}}\hspace*{2mm} \Lambda{\text{\rm -reduced}}
\}
$$
is of cardinal $|\mathcal{A}|\geqslant |\mathcal{S}|-K^{\mathcal{N}(\mathcal{S})}$ which is greater than zero by the above assumption. 

Consider an arbitrary word $W\in \mathcal{A}$. 
By Lemma~\ref{conj_geodes_doubly_Lambda-reduced}, 
there exists a proper terminal subword $T_W$ of $W$ and a word $W_{\bullet}\in (\mathcal{X}\sqcup \mathcal{H})^{\ast}$ such that
\begin{equation}\label{lem_5_2_proof}
\overline{W_{\bullet}}=\overline{T_WWT_W^{-1}}\hspace*{3mm}\text{and}\hspace*{3mm} ||W_{\bullet}||<||W||\leqslant \mathcal{N}(\mathcal{S}),
\end{equation}
$$\overline{T_W}\in \langle X\cup \Omega\rangle\hspace*{3mm}{\text{\rm and}}\hspace*{3mm} 
|\overline{T_W}|_{X\cup \Omega}\leqslant (2M\mathcal{C}+1)||W||\leqslant (2M\mathcal{C}+1)\mathcal{N}(\mathcal{S}).
$$

\medskip

Therefore, there exist at most $K^{\mathcal{N}(\mathcal{S})}$ variants for $\overline{T_W}$ when $W$ runs over $\mathcal{A}$.
Thus, there exists a subset $\mathcal{B}\subseteq \mathcal{A}$ of cardinal $$
|\mathcal{B}|\geqslant |\mathcal{A}|/K^{\mathcal{N}(\mathcal{S})}
$$ 
such that for all words $W\in \mathcal{B}$ the terminal subwords $T_W$ represent the same element of $G$.
We choose some $W_0\in \mathcal{B}$ and set
$$
U=T_{W_0}\hspace*{3mm}
\text{and}\hspace*{3mm}
\mathcal{S}_1=\{W_{\bullet}\,|\, W\in \mathcal{B}\},
$$ 
where $W_{\bullet}$ is associated with $W$ as above.
We show that these $U$ and $\mathcal{S}_1$ satisfy the conclusion of the lemma.
First,
for any $W\in \mathcal{B}$, we have
\begin{equation}\label{about_U}
\overline{U}=\overline{T_{W_0}}=\overline{T_W} \hspace*{3mm}{\text{\rm and}}\hspace*{3mm} ||U||=||T_{W_0}||< ||W_0||\leqslant \mathcal{N}(\mathcal{S}).
\end{equation}

Taking into account the first equation in~\eqref{about_U}, we deduce from~\eqref{lem_5_2_proof}  that
\begin{equation}
\label{5_12_2024}
\overline{\mathcal{S}_1}=\overline{U\mathcal{B}U^{-1}}\subseteq \overline{U\mathcal{S}U^{-1}}\hspace*{3mm}
{\text{and}}\hspace*{3mm} \mathcal{N}(\mathcal{S}_1)<\mathcal{N}(\mathcal{S}). 
\end{equation}

It remains to prove the last inequality in the lemma.
Note that the first equality in~\eqref{5_12_2024} implies that $|\overline{\mathcal{S}_1}|=|\overline{\mathcal{B}}|$. Recall that, by assumption, the map
$\mathcal{S}\rightarrow H$, $W\mapsto \overline{W}$, is injective. Then $|\mathcal{B}|=|\overline{\mathcal{B}}|$, and hence $|\mathcal{S}_1|\geqslant |\overline{\mathcal{S}_1}|=|\mathcal{B}|$.  
The following estimation completes the proof:
$$
|\mathcal{S}_1|\geqslant |\mathcal{B}|\geqslant \frac{1}{K^{\mathcal{N}(\mathcal{S})}}\,|\mathcal{A}|\geqslant \frac{1}{K^{\mathcal{N}(\mathcal{S})}}(|\mathcal{S}|-K^{\mathcal{N}(\mathcal{S})})=\frac{1}{K^{\mathcal{N}(\mathcal{S})}} |\mathcal{S}|-1.
$$
\hfill $\Box$

\medskip







\begin{lemma}\label{not strongly Lambda-reduced_c}
Let $G$ and $\mathcal{P}$ be as above, and let $H$ be a finite subgroup of $G$.
Suppose that $\mathcal{S}$ is a nonempty subset of $(\mathcal{X}\sqcup \mathcal{H})^{\ast}$ that satisfies the following conditions.

\begin{enumerate}
\item[\rm{(a)}] $\mathcal{S}$ bijectively represents a subset of $H$.

\item[\rm{(b)}] Every word $W\in \mathcal{S}$ is geodesic,
has length $||W||\geqslant 3$, and $W_{-},W_{+}\in H_i\setminus \{1\}$ for some $i\in\{1,\dots,m\}$
that does not depend on $W$.
\end{enumerate}

\medskip

\noindent
Then there exist a letter $U\in \mathcal{H}$ and a subset $\mathcal{S}_1\subseteq (\mathcal{X}\sqcup \mathcal{H})^{\ast}$
such that
$$
\overline{\mathcal{S}_1}\subseteq \overline{U\mathcal{S}U^{-1}},\hspace*{5mm}
\mathcal{N}(\mathcal{S}_1)<\mathcal{N}(\mathcal{S}),
\hspace*{5mm}
|\mathcal{S}_1|\geqslant \frac{1}{4K^{2\mathcal{N}(\mathcal{S})}}\,\,
|\mathcal{S}|.
$$
\end{lemma}

\medskip

{\it Proof.} 
By assumption, any word $W\in \mathcal{S}$ can be written in the form $W=W_{-}W'W_{+}$, where $W_{-},W_{+}\in H_i$. We choose an arbitrary word $V\in \mathcal{S}$ 
and write $\mathcal{S}=\mathcal{A}_1\cup \mathcal{A}_2\cup \mathcal{A}_3$, where

$$
\begin{array}{ll}
\mathcal{A}_1= &\{W\in \mathcal{S}\, |\, V_{+}W_{-}=1\} 
,\vspace*{2mm}\\\mathcal{A}_2= & \{W\in \mathcal{S}\,|\, W_{+}V_{-}=1\},\vspace*{2mm}\\
\mathcal{A}_3= & \{W\in \mathcal{S}\, |\, V_{+}W_{-}\neq 1 
\hspace*{2mm}{\text{\rm and }}\hspace*{2mm} W_{+}V_{-}\neq 1\}.
\end{array}
$$ 

\medskip

Then one of the following three cases occurs.

\medskip

{\it Case 1.} Suppose that $|\mathcal{A}_1|\geqslant |\mathcal{S}|/4$. 

We check that the conclusion of lemma is satisfied for
$$
U=V_{+}\hspace*{3mm}
\text{and}\hspace*{3mm}
\mathcal{S}_1=\{(V_{+}WV_{+}^{-1})_{\text {\rm red}}|\, W\in \mathcal{A}_1\}.
$$
The first formula in the conclusion of the lemma follows from
$$
\overline{\mathcal{S}_1}=
\overline{U\mathcal{A}_1U^{-1}}\subseteq \overline{U\mathcal{S}U^{-1}}.
$$
Preparing to prove the second formula, we recall that $V_{+}W_{-}=1$ for $W\in \mathcal{A}_1$. Then 
$$
(V_{+}WV_{+}^{-1})_{\text {\rm red}}=
(V_{+}W_{-}W'W_{+}V_{+}^{-1})_{\text {\rm red}}=
W'(W_{+}V_{+}^{-1})_{\text {\rm red}}.
$$
Therefore
$$||(V_{+}WV_{+}^{-1})_{\text {\rm red}}||\leqslant||W'||+1=||W||-1,
$$
and hence
$$
\mathcal{N}(\mathcal{S}_1)\leqslant \mathcal{N}(\mathcal{A}_1)-1
\leqslant \mathcal{N}(\mathcal{S})-1.
$$

It remains to prove the last inequality of the lemma. To do this, first note that
the map $\mathcal{A}_1\rightarrow \mathcal{S}_1$, $W\mapsto (V_{+}WV_{+}^{-1})_{\rm red}$, is injective.
This follows from the assumption that each $W\in \mathcal{S}$ is geodesic and the fact that geodesic words are reduced. Obviously, this map is surjective. Therefore,
$$
|\mathcal{S}_1|=|\mathcal{A}_1|\geqslant |\mathcal{S}|/4,
$$
and we are done.

\medskip

{\it Case 2.} Suppose that $|\mathcal{A}_2|\geqslant |\mathcal{S}|/4$. 

Then we can set $U=(V_{-})^{-1}$ and 
$
\mathcal{S}_1=\{(V_{-}^{-1}WV_{-})_{\text {\rm red}}|\, W\in \mathcal{A}_2\}.
$
The proof is similar to that in Case~1.

\medskip

{\it Case 3.} Suppose that $|\mathcal{A}_3|\geqslant |\mathcal{S}|/2$.

Let $W$ be an arbitrary word from $\mathcal{A}_3$.
We write 
$$
W=W_{-}w_1\dots w_nW_{+}\hspace*{2mm}{\text {\rm and}}\hspace*{2mm} V=V_{-}u_1\dots u_kV_{+}.
$$
Denote $P=V_{+}W_{-}$ and $Q=W_{+}V_{-}$. Since $W\in \mathcal{A}_3$, we have $P,Q\in H_i\setminus \{1\}$. 
Let $W_{\ast}$ be the result of reducing the word $V_{+}(WV)(V_{+})^{-1}$, i.e., let 
$$
W_{\ast}=Pw_1\dots w_nQu_1\dots u_k.
$$
Note that the word $W_{\ast}$ is cyclically reduced. 
The element $\overline{W_{\ast}}$ has finite order since it is conjugate to the element $\overline{W}\,\overline{V}$ of $H$. 
We subdivide $\mathcal{A}_3$ into two subsets: $\mathcal{A}_3=\mathcal{A}_3'\cup \mathcal{A}_3''$, where
$$
\begin{array}{ll}
\mathcal{A}_3'= & \{W\in \mathcal{A}_3\,|\, W_{\ast}W_{\ast}\hspace*{2mm} {\text{\rm is}}\hspace*{2mm} \Lambda{\text{\rm -reduced}}\},\vspace*{2mm}\\
\mathcal{A}_3''= & \{W\in \mathcal{A}_3\,|\, W_{\ast}W_{\ast}\hspace*{2mm} {\text{\rm is not}}\hspace*{2mm} \Lambda{\text{\rm -reduced}}\}.
\end{array}
$$

\medskip

Then one of the following two cases occurs.

\medskip

{\it Case 3.1.} Suppose that $|\mathcal{A}_3'|\geqslant |\mathcal{A}_3|/2$.

Consider an arbitrary word $W\in \mathcal{A}_3'$. By Lemma~\ref{estimate_Osin}, applied to $W_{\ast}$, we have
$Q\in \langle \Omega\rangle$
and
\begin{equation}\label{Case 3.1}
|Q|_{\Omega} 
\leqslant 2M\mathcal{C}\,||W_{\ast}||=
2M\mathcal{C}\, (||W||+||V||-2)\leqslant 4M\mathcal{C}\cdot \mathcal{N}(\mathcal{S}).
\end{equation}
Hence, when $W$ runs over $\mathcal{A}_3'$, the number of possible $Q$ is at most $K^{2\mathcal{N}(\mathcal{S})}$.
Since $Q=W_{+}V_{-}$ and since $V$ is fixed, the number of possible $W_{+}$ is at most $K^{2\mathcal{N}(\mathcal{S})}$ too. 
Thus, there exists a subset $\mathcal{B}\subseteq \mathcal{A}_3'$ of cardinal 
$$
|\mathcal{B}|\geqslant \frac{1}{K^{2\mathcal{N}(\mathcal{S})}}\, |\mathcal{A}_3'|
$$ 
such that all words $W\in \mathcal{B}$ have the same terminal letter $W_{+}$, say $U$.
We check that the letter $U$ and the set 
$$
\begin{array}{ll}
\mathcal{S}_1 & =\{(W_{+}WW_{+}^{-1})_{\rm red}\,|\, W\in\mathcal{B}\}\vspace*{2mm}\\
& =
\{(UWU^{-1})_{\rm red}\,|\, W\in\mathcal{B}\}
\end{array}
$$
satisfy the conclusion of the lemma. 
As above, the first formula in the conclusion follows from
$$
\overline{\mathcal{S}_1}=
\overline{U\mathcal{B}U^{-1}}\subseteq \overline{U\mathcal{S}U^{-1}}.
$$
Preparing to prove the second formula, we note that 
$$
(W_{+}WW_{+}^{-1})_{\text {\rm red}}=
(W_{+}W_{-}W'W_{+}W_{+}^{-1})_{\text {\rm red}}=
(W_{+}W_{-})_{\text {\rm red}}W'.
$$
Therefore
$$||(W_{+}WW_{+}^{-1})_{\text {\rm red}}||\leqslant||W'||+1=||W||-1,
$$
and hence
$$
\mathcal{N}(\mathcal{S}_1)\leqslant \mathcal{N}(\mathcal{B})-1
\leqslant \mathcal{N}(\mathcal{S})-1.
$$
Finally we prove the last inequality of the lemma:
$$
|\mathcal{S}_1|=|\mathcal{B}|\geqslant \frac{1}{K^{2\mathcal{N}(\mathcal{S})}}\, |\mathcal{A}_3'|\geqslant \frac{1}{2K^{2\mathcal{N}(\mathcal{S})}}\, |\mathcal{A}_3|\geqslant
\frac{1}{4K^{2\mathcal{N}(\mathcal{S})}}\, |\mathcal{S}|.
$$ 

\medskip

{\it Case 3.2.} Suppose that $ |\mathcal{A}_3''|\geqslant |\mathcal{A}_3|/2$.

Consider an arbitrary word $W\in \mathcal{A}_3''$. 
Using above notations, we have  
$$
(W_{\ast})^2=Pw_1\dots w_nQu_1\dots u_kPw_1\dots w_nQu_1\dots u_k.
$$
Let
$L=e_1\dots e_r$ be any path in $\Gamma(G,X\sqcup \mathcal{H})$ with the label $(W_{\ast})^2$, see Fig. 3. 

\medskip

\vspace*{2mm}
\hspace*{-1mm}
\includegraphics[scale=1.0]{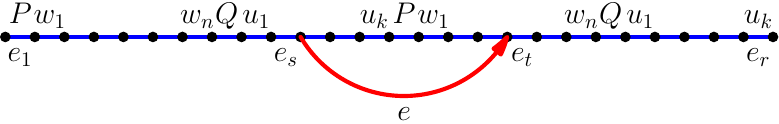}

\vspace*{3mm}

\begin{center}
Fig. 3. Illustration to Case 3.2.
\end{center}

\medskip

Since $W\in \mathcal{A}_3''$, the word $(W_{\ast})^2$ is not $\Lambda$-reduced. Then $L$ contains two connected $H_{\lambda}$-components $e_s$ and $e_t$, $s<t$. Let $e$ be a path in $\Gamma(G,X\sqcup \mathcal{H})$ of length at most~$1$ such that $(e_s)_{+}=e_{-}$,  $(e_t)_{-}=e_{+}$, and $\mathbf{Lab}_G(e)\in H_{\lambda}$. 
Passing to another pair of connected components if necessary, we may assume that all components of the cycle $O=e_{s+1}\dots e_{t-1}\overline{e}$ are isolated.
Since the words $W$ and $V$ are geodesic, the cycle $O$ contains an edge labelled by $P$ or $Q$. The length of $O$ can be estimated as follows:
$$
\ell(O)<||W_{\ast}^2||\leqslant 2(||W||+||V||-2)< 4\mathcal{N}(\mathcal{S}).
$$
Then, by Lemma~\ref{omega_length}, 
we have either $P\in \langle \Omega\rangle$ and $
|P|_{\Omega}\leqslant 4M\mathcal{C}\cdot \mathcal{N}(\mathcal{S}) 
$, or $Q\in \langle \Omega\rangle$ and $
|Q|_{\Omega}\leqslant 4M\mathcal{C}\cdot \mathcal{N}(\mathcal{S}) 
$.
In the second case, we can proceed as in Case 3.1 starting from~\eqref{Case 3.1}. The first case is similar. 
\hfill $\Box$

\section{Proof of the main theorem and corollary}

\subsection{Assumptions that do not change generality.} Recall that $G$ is a relatively hyperbolic group with respect to a collection of subgroups $\{H_{\lambda}\}_{\lambda\in \Lambda}$ and that $\mathcal{P}=\langle X\sqcup\mathcal{H}\,|\, \mathcal{R}\sqcup \mathcal{Q}\rangle$
is a finite relative presentation of $G$ with respect to this collection, where $X$ is a relative generating set. In what follows, we use the notation introduced in Theorem~\ref{finite_subgroups_rel_hyperb_groups}.  
We set 
$$
\Lambda_0=\{\lambda\in \Lambda\,|\, {\text{\rm there exists $R\in \mathcal{R}$ that contains an  $H_{\lambda}$-letter}}\}.
$$
Note that $\Lambda_0$ is finite.

We claim that, without loss of generality, we may assume that 

\begin{enumerate}
\item[{\rm (a)}] $\Lambda=\Lambda_0$ and $\Lambda\neq \emptyset$\\ 
(so we may assume that $\Lambda=\{1,\dots,m\}$ for some natural number $m\geqslant 1$), 

\vspace*{2mm}

\item[{\rm (b)}] $M\geqslant 2$, where $M=\underset{R\in \mathcal{R}}{\max}\, ||R||$.
\end{enumerate}


To show (a), we first note that 
$$
G=G_0\ast (\underset{\lambda\in (\Lambda\setminus \Lambda_0)} {\ast} H_{\lambda}),
$$
where $G_0$ is generated by $X\cup\underset{\lambda\in \Lambda_0}{\cup} H_{\lambda}$. 
The group $G_0$ is relatively hyperbolic with respect to $\{H_{\lambda}\}_{\lambda\in \Lambda_0}$. Furthermore, any finite non-parabolic subgroup of $G$ is
that of $G_0$, and for the induced relative presentation of $G_0$ we can take the same constants $\mathcal{C}\geqslant 1$, 
$\delta\geqslant 1$, $M$, and the set $\Omega$ as for the given relative presentation $\mathcal{P}$ of $G$. 

Therefore, without loss of generality, we may assume that $\Lambda=\Lambda_0$. In particular, $\Lambda$ is finite.
We may also assume that $\Lambda$ is nonempty, otherwise $G$ is hyperbolic and, therefore, any finite subgroup of $G$ has order at most $(2|X|)^{4\delta+2}$ by~\cite{Bog_Ger} (see also~\cite{BH}), and we are done.
Thus, we may assume that condition (a) is satisfied.

(b): If $M\leqslant 1$, then $G$ is the free product of a free group and $\underset{\lambda\in \Lambda}{\ast}H_{\lambda}$. In this case $G$ does not have finite non-peripheral subgroups.  




\subsection{Main lemma.}
In the following definition, we introduce {\it $i$-special subsets} of $G$. This notion will be used in Lemma~\ref{lemma_special_subsets}, the main lemma in the proof.

\begin{definition}
Let $\mathcal{W}_{geod}$
be the set of all geodesic words in the monoid $(\mathcal{X}\sqcup \mathcal{H})^{\ast}$.
In formulations of Lemmas~\ref{strongly Lambda-reduced}
and~\ref{not strongly Lambda-reduced_c}, we used the following two groups of subsets of $\mathcal{W}_{geod}$, respectively:

\begin{enumerate}
\item[$\bullet$] The set $\mathcal{W}_{0}$ consists of all cyclically reduced words of $\mathcal{W}_{geod}$.

\item[$\bullet$] For $i=1,\dots ,m$, the set $\mathcal{W}_i$ consists of all words $W\in \mathcal{W}_{geod}$ satisfying $||W||\geqslant 3$ and
$W_{-}, W_{+}\in H_i$. 
\end{enumerate}

\medskip

Let $i\in \{0,1,\dots ,m\}$.
A subset $A$ of $G$
is called {\it $i$-special} if every element of $A$
can be represented by a word from $\mathcal{W}_i$.
\end{definition}

The following lemma is a direct consequence of Lemmas~\ref{strongly Lambda-reduced} and ~\ref{not strongly Lambda-reduced_c}.
In the following, we write $H^g$ for $gHg^{-1}$.

\begin{lemma}\label{lemma_special_subsets} Let $G$ and $\mathcal{P}$ be as above, and let $H$ be a finite subgroup of $G$.
Suppose that $A\subseteq H$ is a nonempty $i$-special subset for some $i\in \{0,1,\dots, m\}$ and that $$
|A|\geqslant 2K^{2\mathcal{N}(A)}.
$$  
Then there exist an element $g\in G$ of length $|g|_{X\cup \mathcal{H}}< \mathcal{N}(A)$ and a subset $A_1\subseteq H^g$ such that
$$
|A_1|\geqslant \frac{1}{4K^{2\mathcal{N}(A)}}\,|A|\hspace*{4mm}{\text {\rm and}}\hspace*{4mm}
\mathcal{N}(A_1)< 
\mathcal{N}(A).
$$ 
\end{lemma}

\medskip

The following lemma is a direct consequence of the fact that $G\setminus \overset{m}{\underset{j=1}{\cup}} H_j$ is covered by the maximal $i$-special subsets of $G$, $i=0,1,\dots,m$. 

\begin{lemma}\label{2m+2}
Every finite subset $S\subseteq G$ contains a subset $A$ of cardinal at least $|S|/(2m+1)$ such that $A$ is either $i$-special for some $i\in \{0,1,\dots ,m\}$ or lies in $H_j$ for some $j\in\{1,\dots, m\}$. 
\end{lemma}

\subsection{Proof of Theorem~\ref{finite_subgroups_rel_hyperb_groups}}  
In Section 4 we defined the constant
$$
K=(2|X\cup\Omega|)^{2M\mathcal{C}+1}.
$$
In Subsection 5.1, we made assumptions (a) and (b) that do not change generality. In particular, we may assume that $\Lambda=\{1,\dots,m\}$ for some natural number $m\geqslant 1$. 
By (a), we have $|X\cup\Omega|\geqslant |\Omega|\geqslant 
|\Lambda_0|=|\Lambda|=m$. 
Using assumption (b) that $M\geqslant 2$ and condition $\mathcal{C}\geqslant 1$ from Theorem~\ref{finite_subgroups_rel_hyperb_groups}, we deduce that $K\geqslant (2m)^5$. This implies
\begin{equation}\label{K}
K\geqslant 8m+4.
\end{equation}

Let $H$ be a finite non-parabolic subgroup of $G$. Denote $\ell=\mathcal{N}(H)$.
Clearly $\ell\geqslant 1$. 
We will deduce the theorem from the following two claims.

{\it Claim 1.} There exist an element $g\in G$ of length $|g|_{X\cup\mathcal{H}}\leqslant \ell(\ell-1)/2$ and a number $j\in \{1,\dots,m\}$ such that 
\begin{equation}\label{index}
|H^g:(H^g\cap H_j)|\leqslant K^{\mathcal{\ell}^2+2\ell}.
\end{equation}

{\it Proof.}  Consider two cases.

{\it Case 1.} Suppose that $|H|\leqslant 
\overset{\ell}{\underset{i=1}{\prod}}((2m+1)
4K^{2i})$. Then, using~\eqref{K},  we deduce
$$
|H|\leqslant \overset{\ell}{\underset{i=1}{\prod}}
K^{2i+1}=K^{\ell^2+2\ell},
$$
and we are done with $g=1$ and $j=1$.

\medskip

{\it Case 2.} Suppose that $|H|> 
\overset{\ell}{\underset{i=1}{\prod}}((2m+1)
4K^{2i})$.

Then, starting with $S=H$ and applying Lemmas~\ref{2m+2} and~\ref{lemma_special_subsets} alternately several times (at most $\ell$ for each of these lemmas), we conclude that there exists an element $g\in G$ of length $|g|_{X\cup \mathcal{H}}\leqslant \ell(\ell-1)/2$
and a subset $A_{\ast}\subseteq H^g$ such that $A_{\ast}\subseteq H_j$ for some $j\in \{1,\dots ,m\}$ and
$$
|A_{\ast}|\geqslant \frac{1}{
\overset{\ell}{\underset{i=1}{\prod}}((2m+1)
4K^{2i})}\,|H|.
$$
Using~\eqref{K}, we deduce 
$$
|H|/|A_{\ast}|\leqslant  K^{\ell^2+2\ell}.
$$
Then the inequality~\eqref{index} follows from $A_{\ast}\subseteq H^g\cap  H_j$.
\hfill $\Box$

\medskip

By Claim 1 and by the lemma of Poincar\'e, there exists a normal subgroup $N$ of $H^g$ such that 
$N\leqslant H^g\cap H_j$ and 
\begin{equation}\label{index_N}
|H^g:N|\leqslant (K^{\ell^2+2\ell})! 
\end{equation}

{\it Claim 2.} We have
\begin{equation}\label{order_N}
|N|\leqslant K^{(\ell^2+1)}+1.
\end{equation}

{\it Proof}. We may assume that $N$ is nontrivial, otherwise~\eqref{order_N} is evident. 
Since $H$ is non-parabolic, we can choose an element $a\in H^g\setminus H_j$. We have
\begin{equation}\label{length_H^g}
|a|_{X\cup \mathcal{H}}\leqslant \mathcal{N}(H^g)\leqslant 2|g|_{X\cup \mathcal{H}}+\mathcal{N}(H)\leqslant \ell (\ell-1)+\ell=\ell^2.
\end{equation}

Let $b$ be an arbitrary nontrivial element of $N$. We have $a^{-1}ba=b_1$ for some $b_1\in N$. Let $a_1\dots a_k$ be a geodesic word from $(X\sqcup \mathcal{H})^{\ast}$ representing $a$. By~\eqref{length_H^g}, we have 
\begin{equation}\label{k_and_l}
k\leqslant \ell^2.
\end{equation}
Let $P=e_1\dots e_{2k+2}$ be any closed path in
$\Gamma$
with the label $a_k^{-1}\dots a_1^{-1}ba_1\dots a_kb_1^{-1}$, see Fig. 4. In particular, the label of $e_{k+1}$ is $b$.

\medskip

\vspace*{2mm}
\hspace*{26mm}
\includegraphics[scale=1.0]{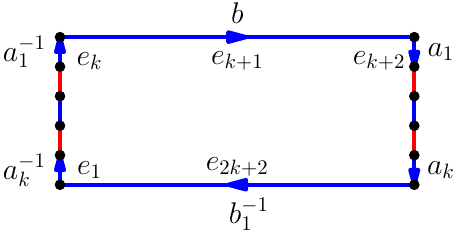}

\vspace*{3mm}

\begin{center}
Fig. 4. The path $P$.
\end{center}

Let $f$ be the $H_j$-component of $P$ containing $e_{k+1}$. Note that $f=e_{k+1}$ if $a_1\notin H_j$ and 
$f=e_ke_{k+1}e_{k+2}$ if $a_1\in H_j$.
Then 
$$
{\text{\bf Lab}}(f)=
\begin{cases}
b & {\text{\rm if}}\hspace*{2mm} a_1\notin H_j,\\
a_1^{-1}ba_1 & {\text{\rm if}}\hspace*{2mm} a_1\in H_j.
\end{cases}
$$
Moreover, $f$ is an isolated $H_j$-component of $P$ because $a\notin H_j$ and the word $a_1\dots a_k$ is geodesic. By Lemma~\ref{omega_length} and~\eqref{k_and_l}, we have
$$
|{\text{\bf Lab}_G}(f)|_{\Omega_j}\leqslant M\mathcal{C}(2k+2)\leqslant 2M\mathcal{C}(\ell^2+1).
$$
Then the number of possible labels of $f$ is at most $(2|\Omega_j|)^{2M\mathcal{C}(\ell^2+1)}$. Since $a_1$ is fixed, the number of possible $b$ is bounded from above by the same number. Then 
$$
|N|\leqslant (2|\Omega_j|)^{2M\mathcal{C}(\ell^2+1)}+1\leqslant K^{(\ell^2+1)}+1. 
$$
\hfill $\Box$ 

\medskip

It follows from~\eqref{K}, \eqref{index_N}, and~\eqref{order_N} that 
$$
\begin{array}{ll}
|H| & \leqslant  (K^{\ell^2+2\ell})! 
(K^{(\ell^2+1)}+1)
\vspace*{2mm}\\
& \leqslant (K^{\ell^2+2\ell})!(K^{\ell^2+2\ell}+1)(K^{\ell^2+2\ell}+2)
\vspace*{2mm}\\
& \leqslant
(K^{(\ell+1)^2})!
\end{array}
$$

Replacing $H$ with a suitable conjugate, we may assume from the beginning that $\ell\leqslant 4\delta+1$, see Theorem~\ref{Bog_Ger} with $Y=X\sqcup \mathcal{H}$. Thus, 
$$
|H|\leqslant(K^{(4\delta+2)^2})!.
$$
The computability of the upper bound under Assumption~\ref{assumprion_rel_hyp} follows from the statements in Subsection~2.6. The last statement of the theorem follows from Claim 1, since if $H_j$ is torsion-free, then $H^g\cap H_j=1$ in~\eqref{index}.  
\hfill $\Box$

\medskip

{\it Proof of Corollary~\ref{corollary}.}
The word problem in $G$ is solvable, since it is solvable in each of the subgroups $H_1,\dots, H_m$ (see~\cite{Farb} or more general Theorem 5.1 in~\cite{Osin_0}).  
 
Let $W$ be a word in the alphabet $Z\sqcup Z^{-1}$ representing a non-parabolic element in~$G$.
We make the list of words $W, W^2,\dots , W^k$, where $k$ is the upper bound from Theorem~\ref{finite_subgroups_rel_hyperb_groups}. If none of these words represents $1$ in $G$, then $\overline{W}$ is the element of infinite order. If some $W^i$ with $i\in \{1,\dots,k\}$ represents $1$, then the smallest such $i$ is the order of $\overline{W}$.

\def\refname{REFERENCES}

\end{document}